\documentclass[runningheads,a4paper]{llncs}

\usepackage{amssymb}
\usepackage{graphicx}

\usepackage{url}
\urldef{\mailsa}\path|d.andreagiovanni@zib.de|
\newcommand{\keywords}[1]{\par\addvspace\baselineskip
\noindent\keywordname\enspace\ignorespaces#1}

\begin{document}

\mainmatter  

\title{On Improving the Capacity of Solving Large-scale Wireless Network Design Problems by Genetic Algorithms\thanks{
This is the authors' final version of the paper published in Di Chio C. et al. (Eds): EvoApplications 2011, LNCS 6625, pp. 11-20, 2011.
DOI: 10.1007/978-3-642-20520-0\_2 $ $.
The final publication is available at Springer via http://dx.doi.org/10.1007/978-3-642-20520-0\_2 $ $.
This work was supported by the \emph{German Federal Ministry of Education and Research} (BMBF), project \emph{ROBUKOM}, grant 03MS616E.}
}

\titlerunning{Improving Solvability of WND Problems by Genetic Algorithms}

%
%
\author{Fabio D'Andreagiovanni%
}
\authorrunning{F. D'Andreagiovanni}

\institute{Konrad-Zuse-Zentrum f\"ur Informationstechnik Berlin (ZIB)\\
Takustrasse 7, 14195 Berlin, Germany\\
\mailsa\\
}

%
%

\toctitle{Lecture Notes in Computer Science}
\tocauthor{Authors' Instructions}
\maketitle

\begin{abstract}
Over the last decade, wireless networks have experienced an impressive growth and now play a main role in many telecommunications systems. As a consequence, scarce radio resources, such as frequencies, became congested and the need for effective and efficient assignment methods arose. In this work, we present a \emph{Genetic Algorithm} for solving large instances of the \emph{Power, Frequency and Modulation Assignment Problem}, arising in the design of wireless networks. To our best knowledge, this is the first Genetic Algorithm that is proposed for such problem. Compared to previous works, our approach allows a wider exploration of the set of power solutions, while eliminating sources of numerical problems. The performance of the algorithm is assessed by tests over a set of large realistic instances of a \emph{Fixed WiMAX Network}.
\keywords{Wireless Network Design, Large-scale Optimization, Genetic Algorithms.}
\end{abstract}

\section{Introduction}

During the last years, wireless communications have experienced an explosive growth thus rapidly leading to a dramatic congestion of radio resources. In such complex scenario, the trial-and-error approach commonly adopted by practitioners to design networks has clearly shown its limitations.
Telecommunications companies and authorities are thus searching for more effective and efficient design approaches, also looking on Optimization (as shown by the recent call for tenders for developing a Digital Video Broadcasting simulator
by the Italian Communications Regulatory Authority \cite{Agcom09}).
Many models and solution methods have been proposed for solving the problem of designing a wireless network. However, solving to optimality the overall problem is still a big challenge in the case of large instances. In this paper we present a Genetic Algorithm for solving the \emph{Power, Frequency and Modulation Assignment Problem}, a relevant problem that arises in the design of wireless networks and captures specific features of Next Generation Networks.

\section{The Wireless Network Design Problem} \label{sec:WND}

For modeling purposes, a wireless network can be described as a set of
transmitters $B$ that provide for a telecommunication service to a set of receivers $T$. A receiver $t \in T$ is said to be {\em covered} or {\em served} if it receives the service within a
minimum level of quality. Transmitters and receivers are
characterized by a location and a number of radio-electrical
parameters (e.g., power emission and frequency).  The {\em Wireless Network Design Problem}  (WND)
consists in establishing the location and suitable values for the parameters of the transmitters with
the goal of maximizing the number of covered receivers or a revenue associated with coverage.

In this work we consider a generalization of the so-called \emph{Power and Frequency Assignment Problem (PFAP)}, a version of the WND that is known to be NP-hard \cite{MaRoSm06}. In addition to power emission and frequency, we also consider the transmission scheme (\emph{burst profile}) as parameter to be established, modeling a feature of Next Generation Networks, such as WiMAX \cite{AnGhMu07,DA10}. We indicate this generalization of the PFAP by the name \emph{Power, Frequency and Modulation Assignment Problem (PFMAP)}.

In the PFMAP two decisions must be taken: (1) establishing the power emission of each transmitter on each available frequency and (2) assigning served receivers to activated transmitters specifying the frequency and the burst profile used to transmit. To model these decisions, as first step we introduce the set $F$ of available frequencies and the set $H$ of available burst profiles.
Each frequency $f \in F$ has a constant bandwidth $D$ and each burst profile $h \in H$ is associated with a \emph{spectral efficiency} $s_h$, which is the bandwidth required to satisfy one unit of traffic demand generated by a receiver.

Then we can introduce two typologies of decision variables, namely:
\begin{itemize}
  \item a continuous power variable $p_b^f \in [0,P_{\max}] \hspace{0.2cm} \forall \hspace{0.1cm} b \in B,  f \in F$ representing the power emission of a transmitter $b$
      on frequency $f$;
  \\
  \item a binary service assignment variable $x_{tb}^{fh} \in \{0,1\} \hspace{0.2cm} \forall \hspace{0.1cm} t \in T, b \in B, f \in F, h \in H$ defined in the following way:
\begin{eqnarray*}
    x_{tb}^{fh} &=& \left\{
                        \begin{array}{lll}
                            1 \hspace{0.3cm} \mbox{ if receiver $t \in T$ is served by transmitter $b \in B$} \\
                            \hspace{0.6cm} \mbox{on frequency $f \in F$ through burst profile $h \in H$} \\
                            0 \hspace{0.3cm} \mbox{ otherwise}
                        \end{array}
                    \right.
\end{eqnarray*}
\end{itemize}

Given a frequency $f \in F$, every receiver $t \in T$ picks out signals from every transmitter $b \in B$ and
the power $P_{b}^{f}(t)$ that $t$ receives from
$b$ on $f$ is proportional to the emitted power $p_b^f$ by a factor
$a_{tb} \in [0,1]$, i.e. $P_{b}^{f}(t) = a_{tb}\cdot p_b^f$.
The factor $a_{tb}$ is called {\em fading coefficient} and
summarizes the reduction in power that a signal experiences while
propagating from $b$ to $t$ \cite{Ra01}.

Among the signals received from transmitters in $B$, a receiver $t$ can select a {\em reference signal} (or {\em server}), which is the
one carrying the service. All the other signals are interfering.
We remark that, since each
transmitter in $B$ is associated with a unique signal,
in what follows we will also refer to $B$ as the set of signals
received by $t$.

A receiver $t$ is regarded as served by the network, specifically
by server $\beta \in B$ on frequency $f \in F$ through burst profile $h \in H$,  if the ratio of the serving power to the
sum of the interfering powers ({\em signal-to-interference ratio}
or \emph{SIR}) is above a threshold $\delta_h$ (\emph{SIR threshold}) whose value depends on the used burst profile $h$ \cite{Ra01}:
\begin{equation}
\label{eq:firstSIRineq} \frac{a_{t \beta} \cdot
p_{\beta}^f}{N+ \sum_{b \in B\setminus\{\beta\}} a_{tb}
\cdot p_b^f}
\hspace{0.1cm}
\geq
\hspace{0.1cm}
\delta_h \;  .
\end{equation}

Note that in the denominator we highlight the presence of the system noise  $N > 0$. By simple algebra operations, inequality (\ref{eq:firstSIRineq}) can be transformed into the
following linear inequality, commonly called \emph{SIR inequality}:
\begin{equation}\label{eq:secondSIRineq}
a_{t \beta} \cdot p_{\beta}^f - \delta_h \sum_{b \in B\setminus\{\beta\}} a_{tb} \cdot p_b^f
\hspace{0.2cm}
\geq
\hspace{0.2cm}
\delta_h \cdot N \;  .
\end{equation}

As we do not know a priori which transmitter $b \in B$ will be the server of a receiver $t \in T$ and which frequency $f \in F$ and burst profile $h \in H$ will be used,
given a receiver $t \in T$ we have one SIR inequality (\ref{eq:secondSIRineq}) for each potential server $\beta \in B$ and potentially usable frequency $f$ and burst profile $h$.
To ensure that $t$ is covered, at least one of such inequalities must be satisfied. This requirement can be equivalently expressed through the following disjunctive constraint:
\begin{equation}\label{eq:disjunctive-SIR}
\bigvee_{\forall (\beta,f,h): \beta \in B, f \in F, h \in H} \left( a_{t \beta} \cdot p_{\beta}^f - \delta_h
\sum_{b \in B\setminus\{\beta\}} a_{tb} \cdot p_b^f
\hspace{0.2cm}
\geq
\hspace{0.2cm}
\delta_h \cdot N \right) \;  .
\end{equation}

Adopting a standard approach used in Mixed-Integer Programming (see \cite{NeWo88}), the above disjunction  can be represented by a family of linear constraints in the $p$ variables by introducing a large positive constant $M$, the so-called \emph{big-M coefficient}. Specifically, given a receiver $t \in T$ we use the assignment variable $x_{t \beta}^{fh}$ to introduce the following constraint for each potential 3-ple $(\beta,f,h)$:
\begin{equation}\label{eq:SIR-BIGM}
 a_{t \beta} \cdot p_{\beta}^f - \delta_h \sum_{b  \in B(t)\setminus\{\beta\}} a_{tb} \cdot p_b^f
\hspace{0.2cm}
+
\hspace{0.2cm}
M  \cdot (1 - x_{t \beta}^{fh})
\hspace{0.2cm}
\geq
\hspace{0.2cm}
\delta_h \cdot N \;  .
\end{equation}

It is easy to check that if
$x_{t \beta}^{fh} = 1$ then (\ref{eq:SIR-BIGM}) reduces to the simple SIR constraint (\ref{eq:secondSIRineq}). If instead $x_{t \beta}^{fh} = 0$ and $M$ is
sufficiently large\footnote{For example, we can set $M = \delta_h \cdot N +
\delta_h \sum_{b \in B\setminus\{\beta\}} a_{tb} \cdot P_{max}$ }, then
(\ref{eq:SIR-BIGM}) is satisfied by any feasible power vector $p$ and
becomes redundant.

By using constraints (\ref{eq:SIR-BIGM}) and by introducing a parameter $r_t$ to denote revenue associated with receiver $t \in T$ (e.g., population, number of customers) , we can define the following natural formulation (BM-PFMAP) for the PFMAP \cite{DA10,DAMaSa10}:
\begin{eqnarray}
    \makebox[13mm][l]{$\max $}
    \makebox[104mm][l]{$
        \displaystyle \sum_{t\in T}
        \sum_{b\in B}
        \sum_{f\in F}
        \sum_{h\in H}
        r_t \cdot x_{t b}^{fh}
        \hspace{4.4cm} (BM-PFMAP)
    $}
    \nonumber
    \\
    \makebox[13mm][l]{s.t.}
    \makebox[104mm][l]{$\displaystyle{a_{t \beta} \cdot p_{\beta}^f -
        \delta_h \sum_{b \in B\setminus\{\beta\}} a_{tb} \cdot p_b^f  + M  \cdot
        (1 - x_{t \beta}^{fh}) \geq \delta_h \cdot N}
    $}
    \nonumber
    \\
    \makebox[45mm][l]{$t\in T,
        b \in B, f \in F, h \in H$}
    \label{eq:BM-SIR}
    \\
    \makebox[59mm][l]{$\displaystyle{\sum_{b\in B} \sum_{f \in F} \sum_{h \in H} x_{tb}^{fh} \leq
    1}$}
    \makebox[45mm][l]{$t \in T$}
    \label{eq:BM-oneserver}
    \\
    \makebox[59mm][l]{$\displaystyle \sum_{t \in T} \sum_{h \in H} d_t \cdot \frac{1}{s_h} \cdot x_{t\beta}^{fh} \leq D$}
    \makebox[45mm][l]{$b \in B, f \in F$}
    \label{eq:BM-capacity}
    \\
    \makebox[59mm][l]{$\displaystyle{p_b^f \in [0, P_{max}]}$}
    \makebox[45mm][l]{$b \in B, f \in F$}
    \label{eq:BM-varP}
    \\
    \makebox[59mm][l]{$x_{t b}^{fh} \in \{0,1\}$}
    \makebox[45mm][l]{$t \in T, b \in B, f \in F, h \in H$}
    \label{eq:BM-varX}
\end{eqnarray}

The objective function is to maximize the total
revenue obtained by serving receivers and constraint (\ref{eq:BM-oneserver}) ensures that each
receiver is served at most once. Each receiver generates a traffic demand $d_t$ and each frequency has a bandwidth equal to $D$. Constraint (\ref{eq:BM-capacity}) ensures that the sum of traffic demands (re-sized by the spectral efficiency $s_h$ of the used burst profile) generated by the receivers served by a transmitter does not exceed the bandwidth of the frequency.
Finally, (\ref{eq:BM-varP}) and (\ref{eq:BM-varX}) define the decision variables of the problem.

\textbf{Drawbacks of the natural formulation.}
The natural formulation (BM-PFMAP) expands a basic model that is widely used for the WND in different application contexts, such as DVB, (e.g., \cite{MaRoSm06}), UMTS (e.g.,
\cite{AmEtAl06,KaKeOl06}) and WiMAX (\cite{DA10,DAMaSa10}).
In principle, such basic model and (BM-PFMAP) can be solved by commercial solvers such as IBM ILOG CPLEX \cite{CPLEX}. However, it is well-known (see \cite{DAMaSa10}) that: (i) the fading coefficients may vary in a wide range leading to (very) ill-conditioned coefficient matrices that make the solution process numerically unstable; (ii) the big-{\em M} coefficients generate poor quality bounds that dramatically reduce the effectiveness of standard solution approach \cite{NeWo88}; (iii) the resulting coverage plans are often unreliable and may contain errors (e.g., \cite{DAMaSa10,KaKeOl06}). In practice, the basic model and (BM-PFMAP) can be solved to optimality only when used for small-sized instances. In the case of large real-life instances, even finding feasible solutions can represent a difficult task, also for state-of-the-art commercial solvers like CPLEX. Though these drawbacks are well-known, it is interesting to note that just a relatively small part of the wide literature devoted to WND has tried to overcome them. We refer the reader to \cite{DA10} for a review of works that have tried to tackle these drawbacks.

\subsection{Contribution of this work and review of related literature} \label{subsect:reviewLit}

In this paper, we develop our original contribution by starting from a recent work, \cite{DAMaSa10}, that proposes a family of strong valid inequalities for tackling the drawbacks of (BM-PFMAP) that we have pointed out. The idea at the basis of \cite{DAMaSa10} is to quit modeling emission power as a continuous variable $p_b$ and to use instead a set of  discrete power levels ${\cal P} = \{P_1, \dots, P_{|{\cal P}|}\}$, with $P_1 = 0$ ({\em switched-off value}), $P_{|{\cal P}|}
= P_{max}$ and $P_i > P_{i-1}$, for $i = 2, \dots, |{\cal P}|$. This basic operation allows the authors to define a family of \emph{lifted GUB cover inequalities} that are used in a solution algorithm that drastically enhances the quality of solutions found.

The solution algorithm proposed in \cite{DAMaSa10} is motivated by a trade-off that arises from discretization: larger sets of discrete levels lead in principle to better solutions, but on the other hand the corresponding 0-1 Linear Program gets larger and harder to solve. The computational experience shows that very good solutions can be found by considering small sets with well-spaced power values, but that no improvement is obtained within the time limit when a number of levels higher than six is considered.

In the present work, we investigate the possibility of using a \emph{Genetic Algorithm (GA)} \cite{Go89} as a fast heuristic to widen the exploration of the discrete power solution space: the aim is to exploit the entire set of discrete power levels and thus to evaluate power configurations with levels not included in the best solutions found in \cite{DAMaSa10}. In particular, our aim is to improve the capacity of solving large realistic instances by finding higher value solutions. We thus design a GA that takes into account the specific features of the PFMAP and we test its performance on the same set of realistic WiMAX instances used in \cite{DAMaSa10}.

Heuristics have been extensively used to tackle large instances of different versions of the WND problem. Two relevant cases are provided by
\cite{AmEtAl06}, where a two-stage Tabu Search algorithm is proposed to solve the base station location and power assignment problem in UMTS networks, and by \cite{MaRoSm06}, where a GRASP algorithm is proposed to solve the PFAP arising in the planning of the Italian National DVB network. The  use of GA to solve versions of the WND is not a novelty as well and many works can be found in literature.
However, to our best knowledge, no GA has been yet developed to solve the PFMAP and the algorithm that we propose is the first for solving this level of generalization of the WND. Until now, GAs were indeed developed to solve just single aspects of the PFMAP: (i) the transmitter location problem (e.g., \cite{ChKiKi08}); (ii) the service assigment problem (e.g., \cite{HuChBa10}); (iii) the frequency assignment problem (e.g., \cite{Co06}); (iv) the power assignment problem (e.g., \cite{SoEtAl02}).
Moreover, we remark that our algorithm is the first to be designed with the specific aim of improving the capacity of solving instances, while tackling the numerical problems pointed out in Section \ref{sec:WND}.
We now proceed to present our original contributions for the WND.

\section{A Genetic Algorithm for the PFMAP} \label{sec:GAalg}

A Genetic Algorithm (GA) is a heuristic method for solving optimization problems that resembles the evolution process of a population of individuals (for a comprehensive introduction to the topic we refer the reader to \cite{Go89}). At any iteration, a GA maintains a population whose individuals represent feasible solutions to the problem. The solution is encoded in a \emph{chromosome} associated with each individual. The \emph{genetic strength} of an individual is evaluated by a fitness function that establishes the quality of the corresponding solution to the problem. A GA starts by defining an initial population, that iteration after iteration changes by crossover, mutation and death of individuals, according to a natural selection Darwinistic mechanism.

We develop a GA for the PFMAP that presents the following general structure:
\begin{enumerate}
  \item Creation of the initial population
  \item UNTIL the arrest condition is not satisfied DO
      \begin{enumerate}
      \item Selection of individuals who generate the offspring
      \item Generation of the offspring by crossover
      \item Mutation of part of the population
      \item Death of part of the population
    \end{enumerate}
\end{enumerate}

We now characterize the elements and the phases presented above for the algorithm (GA-PFMAP) that we propose to solve the PFMAP.

\subsection{Characterization of the population}

\textbf{Individual representation.} As we have explained in Section \ref{subsect:reviewLit}, our aim is to conduct a wider exploration of the power solution space, trying to obtain solutions with higher value. To this end, we establish that the chromosome of an individual corresponds to a power vector $p$ of size $|B|\cdot|F|$. Specifically, the chromosome presents one \emph{locus} for each transmitter $b \in B$ and frequency $f \in F$ and each locus stores the power $p_b^f$ emitted by $b$ on $f$, namely $p = (p_1^1,\ldots,p_1^{|F|}, p_2^1, \ldots, p_2^{|F|}, \ldots, p_{|B|}^{|F|})$. Such power belongs to the set of discrete power levels $\cal{P}$, i.e. $p_b^f \in$ $\cal{P}$ $= \{P_1, \ldots, P_{|\cal{P}|}\}$.

We remark that establishing the power emission $p_b^f \hspace{0.2cm} \forall \hspace{0.1cm} b \in B, \hspace{0.1cm} f \in F$ does not completely characterize a solution of the PFMAP. We indeed have to fix the value of the assignment variables $x_{tb}^{fh}$ and thus we need to set some assignment rule. First, note that given a power vector $p = (p_1^1, p_1^2,\ldots,p_{|B|}^{|F|})$ and a receiver $t \in T$ we can compute the power $P_b^f(t)$ that $t$ receives from $b$ on $f$, $\hspace{0.1cm} \forall b \in B, f \in F$. Through $P_b^f(t)$, if we fix the server $\beta \in B$ of $t$, we can check if there exists a SIR inequality (\ref{eq:secondSIRineq}) that is satisfied for some frequency $f \in F$ and burst profile $h \in H$. We establish the following assignment rule:
as server of $t$ we choose the transmitter $b$ that ensures the highest received power $P_b^f(t)$ on some $f$. This in fact ensures the highest serving power.
Once that the server $\beta$ is chosen, we can identify
the SIR inequalities (\ref{eq:secondSIRineq}) that are satisfied by $p$ when $t$ is served by $\beta$ for some $f \in F$ and $h \in H$. If the SIR inequality is satisfied for a multiplicity of frequencies and/or burst profiles, we first choose as serving frequency $\hat{f}$ the one that ensures the highest value for the left-hand-side of (\ref{eq:secondSIRineq}) and then we choose as burst profile $\hat{h}$ the one that ensures the highest spectral efficiency (see Section \ref{sec:WND}). Thus for $t$ served by $\beta$ we have $x_{t\beta}^{\hat{f}\hat{h}} = 1$ and $x_{t\beta}^{fh} = 0$ $\forall f \in F \setminus \{\hat{f}\}, h \in H \setminus \{\hat{h}\}$.

Note that this last rule may assign a receiver $t$ to a transmitter $\beta$ that violates the capacity constraint (\ref{eq:BM-capacity}) of $\beta$ on the chosen frequency $\hat{f}$. If this is the case, we choose the second best couple of frequency and burst profile according to the rule. If this not possible, the third best and so on. In the end, if there is no capacity left for any valid couple $(f,h)$, $t$ is not considered covered by $\beta$.
\\
\\
\textbf{Fitness function.}
As the aim of the WND is to maximize coverage, we adopt a fitness function that evaluates the coverage ensured by an individual.
Specifically, the fitness $COV(p)$ of an individual is equal to the number of receivers that are covered when the power vector is $p$ and service assignment is done according to the previously introduced rules.
\\
\\
\textbf{Initial population.}
Our aim is to consider all the feasible discrete power levels from the beginning. Therefore, the initial population is represented by the power vectors that are obtained by activating a single transmitter $b \in B$ on a single frequency $f \in F$ at each of the discrete power levels $P_l \in \cal{P}$. For every $b \in B, f \in F$, the corresponding individuals included in the initial population are thus:
$(0,0,\ldots, p_b^f = P_2,\ldots,0,0) \hspace{0.3cm} \cdots \hspace{0.3cm} (0,0,\ldots, p_b^f = P_{|\cal{P}|},\ldots,0,0) \hspace{0.2cm}$.
Note that we exclude the individual corresponding to all transmitters turned off, i.e. $p_b^f = P_1 = 0 \hspace{0.1cm} \forall b \in B, f \in F$. We thus have $|B|\cdot|F|\cdot|L-1|$ initial individuals. We denote the set of individuals representing the population at a generic iteration of the algorithm by $P$.

\subsection{Evolution of the population}

\noindent
\textbf{Selection.}
In order to select the individuals who give birth to the new generation, we adopt a \emph{tournament selection} approach: given the set $P$ of individuals constituting the current population and a value $0 < \alpha < 1$, we first define a number $k \in \mathbb{Z}^+$ of groups including $\lfloor\alpha \cdot |P|\rfloor$ individuals who are randomly extracted from $P$. Then we extract $m < \lfloor\alpha \cdot |P|\rfloor$ individuals with the best fitness from every group. These are the individuals who generate offspring by crossover.
\\
\\
\textbf{Crossover, mutation and death.}
The individuals selected for crossover are randomly paired up to constitute $\lfloor k \cdot  m /2 \rfloor$ couples. Each couple generates two offspring by mixing its chromosome.
Given a couple of parents with power vectors $p1, p2$, the crossover operation consists in mixing power levels that are in the same position of $p1$ and $p2$ to generate two offspring with (possibly) higher fitness power vectors $p3, p4$.

Before presenting the crossover procedure, we need to define a measure that evaluates the impact of crossover on coverage. To this end, let $\Delta\mbox{COV}(p,p_b^f=P_l) \in \mathbb{Z}$ denote the variation in the number of covered receivers caused by changing the power value $p_b^f$ in position $(b,f)$ of vector $p$ to the value $P_l$, while maintaining all the other power values unchanged. We can then propose the following crossover procedure, that concentrates the effort of creating a higher fitness individual on $p3$. At the beginning of the crossover, $p3$ and $p4$ have all elements equal to 0. Then, by following this ordering of indices $(b,f): b \in , f \in F$:
$
(1,1) \mbox{ } (1,2) \mbox{ } \ldots \mbox{ } (1,|F|) \mbox{ } (2,1) \mbox{ } \ldots \mbox{ } (2,|F|) \mbox{ } \ldots \mbox{ } (|B|,1) \mbox{ } \ldots \mbox{ } (|B|,|F|)
$,
each null value inherits the power value in the same position of one of the two parents.

We now present the crossover rule for a generic position $(\beta,\phi)$. For indices $(b,f): b < \beta, f < \phi$, the crossover was already executed and thus the offspring vectors $p3,p4$ present power levels inherited by the parents $p1, p2$. Power levels of $p3,p4$ in positions $(b,f): b \geq \beta, f \geq \phi$ are instead still equal to zero. The rule to establish the power value inherited by $p3,p4$ in $(\beta,\phi)$  is the following:
\begin{eqnarray*}
    p3_\beta^\phi &=& \left\{
                        \begin{array}{lll}
                            p1_\beta^\phi \hspace{0.3cm} \mbox{ if $\Delta\mbox{COV}(p3,p3_\beta^\phi=p1_\beta^\phi) \geq \Delta\mbox{COV}(p3,p3_\beta^\phi=p2_\beta^\phi)$ }
                            \\
                            p2_\beta^\phi \hspace{0.3cm} \mbox{ otherwise}
                        \end{array}
                    \right.
\\
\\
    p4_\beta^\phi &=& \left\{
                        \begin{array}{lll}
                            p1_\beta^\phi \hspace{0.3cm} \mbox{ if $\Delta\mbox{COV}(p3,p3_\beta^\phi=p1_\beta^\phi) < \Delta\mbox{COV}(p3,p3_\beta^\phi=p2_\beta^\phi)$ }
                            \\
                            p2_\beta^\phi \hspace{0.3cm} \mbox{ otherwise}
                        \end{array}
                    \right.
\end{eqnarray*}

This ensures that, at any step of the crossover procedure, offspring $p3$ inherits the power level of the parent that allows the most favourable variation $\Delta COV$ in coverage.

In addition to crossover, we also allow to alter the power vector of single individuals by \emph{mutation}. This introduces new genetic information in the population and helps to widen the solution space exploration and to avoid entrapment in local optima. At any iteration, a number of individuals $\lfloor\gamma \cdot |P|\rfloor$ with $0 < \gamma < 1$ is randomly chosen. Then, still by random selection, $|F|$ power levels corresponding with different frequencies are reduced to the immediately lower power level allowed in $\cal{P}$. This mutation rule is set with the aim of defining new power vectors that have lower powers but ensure the same coverage. The reduction in power is generally desirable as a signal that is useful for a receiver may be interfering for a different receiver.

Finally, after crossover and mutation, the weakest individuals \emph{die} and are removed from $P$. Specifically, we choose to select and remove the $2 \cdot \lfloor k \cdot  m /2 \rfloor$ individuals with the worst fitness function. The size of \emph{P} is thus maintained constant over all the iterations.

\section{Computational experience} \label{sec:comput}

We test the performance of our GA on a set of 15
realistic instances, developed with the Technical Strategy \&
Innovations Unit of British Telecom Italia (BT Italia SpA). All the instances refers to a \emph{Fixed WiMAX Network} \cite{AnGhMu07}, deployable in an urban area corresponding to a
residential neighborhood of Rome (Italy). The instances consider various scenarios with up to $|T| = 529$ receivers, $|B| = 36$ transmitters, $|F| = 3$ frequencies, $|H| = 4$ burst profiles (see Table \ref{tab:comparisons}). This leads to large formulations (BM-PFMAP) that are very hard to solve. For a detailed description of the instances, we refer the reader to \cite{DAMaSa10}.

For each instance, we run the proposed algorithm (GA-PFMAP) 50 times with a time limit of 1 hour by using a machine with a 1.80 GHz Intel Core 2 Duo processor and 2 GB of RAM. Each tournament selection involves $k = 20$ groups that include a fraction $\alpha = 0.05$ of the population $P$. The best $m = 8$ individuals of each group are selected for crossover and, after the generation of the new individuals, mutation affects a fraction $\gamma = 0.1$ of the population.

In Table \ref{tab:comparisons}, we compare the value of the best solution obtained through the three approaches that we consider, namely the direct solution of (BM-PFMAP) by ILOG Cplex 10.1, the solution of the Power-Indexed formulation by the algorithm WPLAN \cite{DAMaSa10} and the solution of (BM-PFMAP) by the proposed algorithm (GA-PFMAP). Results for (BM-PFMAP) and WPLAN  are derived from \cite{DAMaSa10}. The presence of two values in some lines of the column of (BM-PFMAP) indicates that the coverage plans returned by Cplex contain errors and some receivers are actually not covered. We remark that (GA-PFMAP) provides solutions that always ensure a higher coverage than (BM-PFMAP) and without coverage errors. Making a comparison with WPLAN, we instead note that (GA-PFMAP), though in some cases finds solutions with lower coverage, for most of the cases is able to find solutions that ensure an equal or higher number of covered receivers than WPLAN. This is particularly evident for instances that seems to be very hard to solve through (BM-PFMAP). The algorithm is thus effective and is worth of further investigations.

\begin{table}[tbp]
\caption{Comparisons between (BM) and WPLAN formulations}
\label{tab:comparisons}
\centering
\begin{center}
\begin{tabular}{|ccccc||cccc|}
  \hline
  & & & & & & & $|\mbox{T*}|$ &
  \\
  \raisebox{1.5ex}{ ID } & \raisebox{1.5ex}{$\mbox{ } |\mbox{T}|$ \mbox{ }} &
  \raisebox{1.5ex}{$\mbox{ } |\mbox{B}|$ \mbox{ }} &
  \raisebox{1.5ex}{$\mbox{ } |\mbox{F}|$ \mbox{ }} &
  \raisebox{1.5ex}{$\mbox{ } |\mbox{H}|$ \mbox{ }} &
  & (BM-PFMAP) &
  $\mbox{ }$ WPLAN \cite{DAMaSa10} $\mbox{ }$ &
  (GA-PFMAP) $\mbox{ }$
  \\
  \hline
  \hline
  S1 & 100 & 12 & 1 & 1 & $\mbox{ }$ & 63 (78) & 74 & 70
  \\
  S2 & 169 & 12 & 1 & 1 & $\mbox{ }$ & 99 (100) & 107  & 103
  \\
  S3 & 196 & 12 & 1 & 1 &  $\mbox{ }$ & 108 & 113 & 113
  \\
  S4 & 225 & 12 & 1 & 1 & $\mbox{ }$ &  93 & 111  &  115
  \\
  S5 & 289 & 12 & 1 & 1 & $\mbox{ }$ &  77  & 86 & 88
  \\
  S6 & 361 & 12 & 1 & 1 & $\mbox{ }$ & 154  & 170 &  175
  \\
  S7 & 400 & 18 & 1 & 1 & $\mbox{ }$ & 259 (266) & 341 & 319
  \\
  \hline
  R1 & 400 & 18 & 3 & 4 &  $\mbox{ }$ & 370 &   400 & 400
  \\
  R2 & 441 & 18 & 3 & 4 &  $\mbox{ }$ & 302 (303)  &  441 &  441
  \\
  R3 & 484 & 27 & 3 & 4 &  $\mbox{ }$ &  99 (99)  & 427 & 434
  \\
  R4 & 529 & 27 & 3 & 4 &  $\mbox{ }$ & 283 (286)  &  529 & 462
  \\
  \hline
  Q1 & 400 & 36 & 1 & 4 & $\mbox{ }$ &  0  & 67 & 72
  \\
  Q2 & 441 & 36 & 1 & 4 & $\mbox{ }$ &  191  & 211 & 222
  \\
  Q3 & 484 & 36 & 1 & 4 & $\mbox{ }$ &  226  & 463 & 466
  \\
  Q4 & 529 & 36 & 1 & 4 & $\mbox{ }$ & 145 (147) & 491 & 491
  \\
\hline
\end{tabular}
\end{center}
\end{table}

\section{Conclusion and future work}  \label{sec:end}

We presented a Genetic Algorithm (GA) to tackle large realistic instances of a relevant problem arising in wireless network design. We showed that a GA helps to improve the value of solutions found through a wider exploration of the power space. A future research path could be represented by the integration of a refined GA into an exact solution method.
It is indeed common experience that the combination of fast heuristics with Mixed-Integer Linear Programming leads to a great reduction in the running times w.r.t pure exact optimization methods.
A sensitivity analysis of the GA parameters and a study on the impact of different starting conditions and selection strategies would also constitute important subjects of investigations.

\subsubsection*{Acknowledgments.} The author thanks Carlo Mannino and Antonella Nardin for fruitful discussions. Thanks go also to the three anonymous referees who provided valuable comments.

\end{document}